\input amstex
\documentstyle{amsppt}
\magnification=\magstep1 \pagewidth{5.2 in} \pageheight{6.7in}
\abovedisplayskip=10pt \belowdisplayskip=10pt
\parskip=4pt
\parindent=5mm
\baselineskip=2pt \NoBlackBoxes

\topmatter

\title  $q$-analogue of Euler-Barnes multiple zeta functions
\endtitle
\author  Taekyun Kim\endauthor
\affil\rm{ {Institute of Science Education,}\\
{ Kongju National University, Kongju 314-701, S. Korea}\\
{e-mail: tkim$\@$kongju.ac.kr (or tkim64$\@$hanmail.net)}}\\\\
\endaffil

\define\BZ{\Bbb Z}
\define\BQ{\Bbb Q}
\define\BC{\Bbb C}
\define\BN{\Bbb N}
\define\BR{\Bbb R}

\thanks 2000 Mathematics Subject Classification : 11K31, 11SM06 \endthanks
\keywords Gamma function, Riemann zeta-function, $q$-Euler numbers
\endkeywords
 \abstract{  Recently (see [1]) I has introduced an interesting the Euler-Barnes multiple zeta function.
  In this paper we construct the $q$-analogue of
  Euler-Barnes multiple zeta function which interpolates the $q$-analogue of Frobenius-Euler numbers of higher order at negative
  integers. }
\endabstract
\rightheadtext{ T.Kim } \leftheadtext{ $q$-analogue of
Euler-Barnes' zeta functions} \TagsOnRight
\endtopmatter

\document

\head \S 1. Introduction \endhead

\vskip0.4cm

 Throughout this paper $\BZ,\;\BQ,\; \BR,$ and  $\BC$ will denote the ring of rational integers, the field of
 rational numbers, the field of real numbers and the field of complex numbers, respectively.
In this paper, we use the notation

$$[x]_q =  { {1- q^x} \over  {1-q }}  , \;\;\; cf.
[1,2,3,4,5,6].$$

The ordinary Euler numbers $E_m$ are defined by the generating
function in the complex number field as
$$ { 2 \over {e^t +1 }} =e^{Et}= \sum_{m=0}^\infty E_m {{t^m } \over {m!}} ,
\;\; (|t| < \pi) , $$ where we use the usual convention about
replacing $E^m$ by $E_m (m
>0)$, symbolically.

Let $u$ be algebraic in complex number field. Then Frobenius-Euler
numbers are defined by
$$ { {1-u} \over {e^t -u }} = \sum_{n=0}^\infty H_n (u) {{t^n } \over {n!}} ,
\;\; (|t| < \pi) ,\;\; cf. [1]. $$ Note that $H_n (-1) = E_n $.

The Bernoulli numbers $B_n $ are defined as
$$ { t \over {e^t -1 }} = \sum_{n=0}^\infty B_n {{t^n } \over {n!}} ,
\;\; (|t| < \pi) , \;\; cf.[1,2,3,4].$$ Then we note that
$$ H_n (-1) = \sum_{k=0}^n \binom {n+1}k 2^k B_k ,$$
where $\binom nk $ is a binomial coefficient, cf. [1,2,3,4].

Recently I has introduced the Euler-Barnes zeta function in [1].
In this paper we give the q-analogue of Euler-Barnes' zeta
function which interpolates the Frobenius-Euler numbers and
polynomials of higher order. Finally, we give some interesting
formulas for the $q$-analogue of Frobenius-Euler numbers and
polynomials. \vskip0.4cm

\head \S 2. $q$-Euler numbers and polynomials
\endhead

\vskip 0.4cm

For $q \in \BC$ with $|q|<1$, let $u$ be algebraic number in the
complex number field with $|u|<1$. Then we consider the $q$-Euler
numbers which are defined by the generating function in the
complex number field as
$$ F_{u^{-1}, q} (t) = (1-u) \sum_{l=0}^\infty u^l e^{[l]_q t}
=e^{H_q (u)t}= \sum_{n=0}^\infty H_{n,q} (u^{-1}) {{t^n}\over
{n!}},$$ where we use the usual convention about replacing $H_q^n
(u)$ by $H_{n,q}(u),\;(n \geq 0),$ symbolically. By simple
calculation, we note that
$$\aligned
(1-u) \sum_{l=0}^\infty u^l e^{[l]_q t}
&=(1-u) e^{{1 \over {1-q}}t} \sum_{l=0}^\infty u^l e^{-{{q^l}\over {1-q}}t}
\\
&= (1-u) e^{{1 \over {1-q}}t} \sum_{l=0}^\infty u^l
\sum_{j=0}^\infty \left({1 \over {1-q}}\right)^j (-1)^j q^{lj}
{{t^j}\over
{j!}}\\
 & = (1-u) e^{{1 \over {1-q}}t}  \sum_{j=0}^\infty \left({ 1 \over {1-q}} \right)^j (-1)^j { 1 \over {1-u q^j } } {{t^j} \over
{j!}}.
\endaligned  \tag2$$

By (1) and (2), we easily see that
$$ H_{n,q} (u^{-1} )= {{1-u} \over {(1-q)^n }}  \sum_{l=0}^n \binom nl (-1)^l {1\over{1-uq^l}} .
\tag 3$$

Now, we define $q$-Euler polynomials as
$$\aligned
F_{u^{-1} , q} (x,t)&=e^{[x]_q t} F_{u^{-1},q} (q^x t)\\
&= (1-u) \sum_{l=0}^\infty u^l e^{[x+l]_q t}\\
&= \sum_{n=0}^\infty H_{n,q}(u^{-1} ,x) {{t^n } \over {n!}}.
\endaligned  \tag4$$
By (3) and (4), we easily see that
$$\aligned
H_{n,q} (u^{-1} , x) &=  \sum_{l=0}^n \binom nl (-1)^l q^{lx}
\left({1
\over {1-uq^l}} \right) \left( {{1-u}\over {1-q^n}} \right)\\
&= \sum_{l=0}^n \binom nl [x]_q^{n-l} q^{lx} H_{l,q} (u^{-1} ) .
\endaligned  \tag5$$

For $s \in \BC$, let us consider the below complex integral.
$$ {1 \over {1-u}} {1 \over {\Gamma(s)}} \int_0^\infty F_{u^{-1} ,
q}(x, -t) t^{s-1} dt=\sum_{l=0}^\infty u^l {l\over {\Gamma(s)}}
\int_0^\infty e^{-[x+l]_q t} t^{s-1} dt=\sum_{l=0}^\infty
{{u^l}\over {[x+l]_q^s}},  \tag 6$$ where $\Gamma(s)$ is the gamma
function.

Thus we can construct the $q$-analogue of Euler-Hurwitz zeta
function, cf. [1], as
$$ \zeta_q ( u| s,x)= \sum_{n=0}^\infty {{u^n } \over {[n+x]_q^s}}, \text{ for $s\in\Bbb C$} .
\tag 7$$ It is easy to see that $\zeta_q (u| s,x)$ is analytic
continuation for $Re(s)>1$. Let us define $\zeta_q (u| s)$ ,which
is called the $q$-analogue of Euler-Riemann zeta function , as
$$ \zeta_q ( u| s)= \sum_{l=1}^\infty {{u^l } \over {[l]_q^s}}, \text{ for $s\in\Bbb C$} .
\tag 8$$

 For $n \in \BN$, we note that
$$ \zeta_q (u| -n,x)= {1 \over {1-u}} H_{n,q} (u^{-1} , x) =
{{u^{-1}} \over {u^{-1} -1}} H_{n,q} (u^{-1}, x). \tag 9$$ We now
define the $q$-analogue of Barnes-Euler polynomials as
$$ \aligned
 F_{u^{-1},q}^{(r)} (x,t) &= (1-u)^r \sum_{n_1, \cdots, n_r
 =0}^\infty u^{n_1 + \cdots + n_r} e^{[x+ n_1 + \cdots + n_r]_q t
 }\\
 &=\sum_{n=0}^\infty H_{n,q}^{(r)} (u^{-1}, x) {{t^n}\over {n!}} .
 \endaligned \tag 10$$
Thus we may consider the below complex analytic $r$-ple $q$-Euler
zeta function (or the $q$-analogue of Barnes-Euler $r$-ple zeta
function).
$$ \zeta_{r,q}(u|s,x) = \sum_{n_1 , \cdots , n_r =0}^\infty
{{u^{n_1 + \cdots + n_r}} \over {[x+n_1 +\cdots +n_r]_q^s}}, \;
\;\;  for\;\; \Re(x)>0, \text{ $s\in\Bbb C$}.$$ Note that
$\lim_{q\rightarrow 1}\zeta_{r,q}(u|s,x)=\zeta_{r}(s,x, u|1,
\cdots, 1), \text{ see [1]. }$ Let $\zeta_{r,q} (u|s)=u^r
\zeta_{r,q}(u|s,r)$. Analytic continuation and special values of
$\zeta_{r,q} (u|s,x)$ are given by the below complex integral
representation.
$$ \zeta_{r,q} (u|s,x)= \left( {1\over{\Gamma(s)}} \int_0^\infty
F_{u^{-1},q}^{(r)} (x, -t) t^{s-1} dt \right) {1 \over {(1-u)^r}},
\text{ for $s\in\Bbb C$ } .$$  By (10), we easily see that
$$\zeta_{r,q}(u| -n,x) ={1 \over {(1-u)^r}} H_{n,q}^{(r)} (u^{-1},
x), \text{ for $n\in\Bbb N$} .$$ Let $\chi$ be the Dirichlet
character with conductor $d\in\Bbb N .$
 We now define the generalized $q$-Euler numbers attached to $\chi$ as
 $$F_{u,\chi: q}(t)=(1-u)\sum_{n=0}^{\infty}e^{[n]_q t}\chi(n)u^n
 =\sum_{n=0}^{\infty} H_{n,\chi, q}(u^{-1})\frac{t^n}{n!}. \tag11
 $$
Note that
$$H_{n,\chi,q}(u^{-1})=\frac{1-u}{1-u^d}[d]_q^n\sum_{a=0}^{d-1}\chi(a)u^a
H_{n,q^d}(u^{-d},\frac{a}{d}).$$ For $s\in\Bbb C, $ let us define
the $q$-$l$-function as
$$l_q(s, \chi)=\sum_{n=1}^{\infty}\frac{\chi(n)u^n}{[n]_q^s}=\frac{1}{(1-u)\Gamma(s)}
\int_{0}^{\infty}F_{u,\chi:q}(-t)t^{s-1} dt. $$ Then we easily see
that
$$l_q(-n, \chi)=\frac{1}{1-u}H_{n, \chi, q}(u^{-1}), \text{ for $n\in \Bbb N$}.$$

\enddocument